\documentclass[12pt,draft]{extartic}

\usepackage{amsmath,amsthm,amssymb,calc}
\usepackage[dvips]{graphics, graphicx}

\usepackage[english]{babel}

\def\eps{\varepsilon}

\newcounter{num}[section]

\newcommand{\Th}{\refstepcounter{num}
{\bf Theorem \arabic{section}.\arabic{num} }}

\newcommand{\Lemma}{\refstepcounter{num}
{\bf Lemma \arabic{section}.\arabic{num} }}

\newcommand{\Pred}{\refstepcounter{num}
{\bf Proposition \arabic{section}.\arabic{num} }}

\newcommand{\Cor}{\refstepcounter{num}
{\bf Corollary \arabic{section}.\arabic{num} }}

\newcommand{\Note}{\refstepcounter{num}
{\it Note \arabic{section}.\arabic{num} }}

\newcommand{\Exm}{\refstepcounter{num}
{\bf Example \arabic{section}.\arabic{num} }}

\newcommand{\Def}{\refstepcounter{num}
{\it Definition \arabic{section}.\arabic{num} }}

\newcommand{\Proof}{{\bf Proof. }}

\def\eps{\varepsilon}
\def\_phi{\varphi}

\def\a{\alpha}
\def\d{\delta}
\def\la{\lambda}

\def\F{\widehat}

\def\t{\tilde}
\def\o{\omega}
\def\ov{\overline}

\def\C{{\mathbb C}}
\def\R{{\mathbb R}}
\def\E{{\mathbb E}}

\def\Z_N{{\mathbb Z}_N}
\def\Z{{\mathbb Z}}

\def\Gr{{\mathbf G}}

\def\E{\mathsf {E}}

\def\D{{\mathbb D}}

\def\supp{{\rm supp\,}}

\def\D{{\mathbb D}}

\def\supp{{\rm supp\,}}

\def\C{{\mathbb C}}
\def\Gr{{\mathbf G}}

\def\FF{\widehat}
\def\c{\circ}
\def\D{\Delta}
\def\Cf{{\mathcal C}}

\def\F{\mathbb {F}}

\def\no{\noindent}

\author{Shkredov I.D.}
\title{ Sumsets in quadratic residues
\footnote{ This work was supported by grant RFFI NN
11-01-00759, Russian Government project 11.G34.31.0053, Federal
Program "Scientific and scientific--pedagogical staff of
innovative Russia" 2009--2013, grant mol\underline{ }a\underline{
}ved 12--01--33080 and grant Leading Scientific Schools N
2519.2012.1.}}
\date{}
\begin{document}
\maketitle

\begin{center}
 Annotation.
\end{center}

{\it \small
    We describe all sets $A \subseteq \F_p$ which represent the quadratic residues $R \subseteq \F_p$
    as $R=A+A$ and $R=A\hat{+} A$.
    Also, we consider the case of an approximate equality
    $R \approx A+A$ and $R \approx  A\hat{+} A$ and prove that $A$ has a structure in the situation.
}
\\

\refstepcounter{section}
\label{sec:introduction}

\begin{center}
{\bf \large
    \arabic{section}. Introduction}
\end{center}

Let $p$ be a prime number, $\F_p$ be the finite  field,
and $R$ be a the set of all quadratic residues.
In other words, $R$ is the set of all squares in $\F_p \setminus \{ 0 \}$.
There are many interesting conjectures about the set $R$ (see e.g. \cite{Palo_Alto_problems,Karatsuba,Sarkozy_R}).
We begin with a hypothesis  of A. S\'{a}rk\"{o}zy \cite{Sarkozy_R}.

{\bf Conjecture A.}
{\it
    The set $R$ cannot be represented as a sumset $R=A+B$,
    where the cardinality  of each set $A,B$ is at least $2$.
}
\label{conj:Sarkozy}

Here, as usual,
$$
    A+B := \{ a+b ~:~ a\in A,\, b\in B \} \,.
$$
In the paper we will also need in another definition.
For a set $A$ put
$$
    A \hat{+} A := \{ a+a' ~:~ a,a'\in A,\, a\neq a' \} \,.
$$

In \cite{Sarkozy_R} the following result was obtained (see also \cite{Sarkozy_R_g}).

\Th
{\it
    Let $p$ be a prime number.
    Suppose that $R=A+B$, $|A|,|B| \ge 2$. Then
$$
     \frac{p^{1/2}}{3\log p} < |A|,|B| < p^{1/2} \log p \,.
$$
}
\label{t:Sarkozy}

Some generalizations and improvements of Theorem \ref{t:Sarkozy} can be found in \cite{Sarkozy_R_g,Shp_A+B}.

Another well--known conjecture (see e.g. \cite{Palo_Alto_problems}) asks, in particular, the following.

{\bf Conjecture B.}
{\it
    Let $\eps \in (0,1)$ be a real number and $p$ be a sufficiently large prime number.
    Suppose that $A+A \subseteq R$ or $A-A \subseteq R \sqcup \{ 0 \}$.
    Then $|A| \ll p^{\eps}$.
}
\label{conj:A+A}

Some results in the direction can be found in \cite{Ruzsa_sqrt,Paley1,Paley2}.
At the moment the best known bound has the form $|A| \ll \sqrt{p}$, see e.g. \cite{Ruzsa_sqrt}.
A lower bound for the case $A-A \subseteq R \sqcup \{ 0 \}$ is due to S. Graham and C. Ringrose \cite{Paley2}.
It asserts us that $|A| \gg \log p \cdot \log \log \log p$ for infinitely many primes $p$.
A uniform lower bound of the form $|A| \ge \left( \frac{1}{2} + o(1) \right) \log p$,
where $p$ is an arbitrary prime, can be found in \cite{Paley1}.
Conjecture B can be reformulated in terms of the clique number of the Paley graph $P_p$, see e.g. \cite{Bollobas}.

Exponential sums with multiplicative characters over sumsets have been studied by various authors
(see e.g. \cite{Chang_char1},
\cite{DE1}, \cite{DE2}, \cite{Karatsuba_-2}, \cite{Karatsuba_-1}, \cite{Karatsuba}).
The classical general result of \cite{DE1}, \cite{DE2} states
\begin{equation}\label{f:DE}
    \sigma (A,B) := \left| \sum_{x,y} \chi(x+y) B(x) A(y) \right|
        \le
            \sqrt{|A||B|p}
\end{equation}
for any sets $A,B\subseteq \F_p$ and an arbitrary non principal multiplicative character $\chi$.
The bound is nontrivial if $|A||B| > p^{1+\d}$, $\d>0$.
A.A. Karatsuba and M.--C. Chang to contribute significantly in to the theory of such exponential sums.
For example, Karatsuba proved a non--trivial upper bound for the sum $\sigma (A,B)$
provided by $|A|>p^{\eps_1}, |B| > p^{1/2+\eps_2}$, $\eps_1,\eps_2>0$.
Chang obtained plenty results for specific $A$ and $B$, e.g. $A$ has small sumset or, conversely, $A$ is a well--spaced set
(
see \cite{Chang_char1} for details).
In his survey \cite{Karatsuba} Karatsuba formulated a hypothesis (see Problem 6).

{\bf Conjecture C.}
{\it
    Let $A,B\subseteq \F_p$, $|A|,|B| \sim \sqrt{p}$.
    Then
\begin{equation}\label{f:conj_C}
    \left| \sum_{x,y} \chi(x+y) B(x) A(y) \right|
        \le
            c(\d) |A| |B| p^{-\d} \,, \quad \d>0 \,.
\end{equation}
}
\label{conj:Karacuba}

The last hypopiesis in the section which again is called Paley graph conjecture (see e.g. \cite{Chang_char1}) predicts
a non--trivial upper bound for the sum
(\ref{f:DE}) for $|A|,|B|>p^{\eps}$, $\eps>0$.

{\bf Conjecture D.}
{\it
    Let $A,B\subseteq \F_p$, $|A|,|B| > p^{\eps}$, $\eps>0$.
    Then
\begin{equation}\label{f:conj_D}
    \left| \sum_{x,y} \chi(x+y) B(x) A(y) \right|
        \le
            c(\eps) |A| |B| p^{-\d} \,, \quad \d = \d(\eps) >0 \,.
\end{equation}
}
\label{conj:Paley}

Clearly, Conjecture D is the strongest one and implies all another.
Trivially, Conjecture A follows from
Conjectures B or Conjecture C.
It is known that  the correspondent functional version of Conjecture C is false, see e.g. section \ref{sec:appendix}.
In the paper we give partial answer on Conjecture A.
Let us formulate our main result.

\Th
{\it
    Let $p$ be a prime number,  $R \subseteq \F_p$ be the set of quadratic residues and $A\subseteq \F_p$ be a set. \\
${\bf 1)}~$ If $A+A = R$ then $p=3$ and  $A=\{ 2 \}$. \\
${\bf 2)}~$ If $A \hat{+} A = R$ then $p=3,7,13$
and there are just four possibilities for $A$, see Example \ref{exm:pd_sets}. \\
}
\label{t:main_intr}

Also, we improve a little bit Theorem \ref{t:Sarkozy} as well as consider the cases of
approximate equalities, in some sense,
that is
$A+A \approx R$ and $A \hat{+} A \approx R$,  see sections \ref{sec:proof}, \ref{sec:proof2}.
Note that an improvement of S\'{a}rk\"{o}zy's theorem was  obtained independently
by I. Shparlinski in \cite{Shp_A+B} using Karatsuba's bound from \cite{Karatsuba_-1}.
As for Conjecture B, we
reprove
a recent result of C. Bachoc, M. Matolcsi, I.Z. Rusza \cite{Ruzsa_sqrt}
in the direction.
Interestingly, our method does not use the main lemma of
paper \cite{Ruzsa_sqrt}.

In their proof the authors \cite{Sarkozy_R,Sarkozy_R_g,Shp_A+B}
used the well--known Weil bound for exponential sums with
multiplicative characters (see e.g. \cite{Johnsen})
\begin{equation}\label{f:Weil}
    \left| \sum_x \chi(x) \chi(x+x_1) \dots \chi(x+x_d) \right|
        \le
            (d-1) \sqrt{p}
\end{equation}
for different nonzero $x_1,\dots, x_d \in \F_p$ and  any non principal character $\chi$,
 as well as some combinatorial tools.
Our main idea exploit the fact that quadratic residues are "more random then a random set"\,
and, hence, structured.
The last statement is the most transparent
at the case $p\equiv 3 \pmod 4$
(or, equivalently, $\binom{-1}{p} = -1$, in terms of Legendre symbol).
It is known that $R$ is a perfect difference set (see e.g. Lemma \ref{l:R_conv} from section \ref{sec:definitions})
in the situation, that is
the number of solutions $x=a-b$ with $a,b\in R$ does not depend on $x\neq 0$.
Of course a random set of density $1/2$ satisfies the property with probability zero.
So, because of we use such properties of $R$ instead of its random behavior
it is very natural that perfect difference sets appear in our proofs.
For example, all sets $A$ from the second part of Theorem \ref{t:main_intr} turn out to be perfect difference sets.


The author is grateful to S. Yekhanin for very useful discussions
and S. Konyagin for intense interest to  the paper.

\refstepcounter{section}
\label{sec:definitions}

\begin{center}
{\bf \large
    \arabic{section}. Notation and auxiliary results}
\end{center}

We start  with definitions and notation used in the paper.
Let $\Gr$ be a finite abelian group.
It is well--known~\cite{Rudin_book} that the dual group $\FF{\Gr}$ is isomorphic to $\Gr$ in the case.
Let $f$ be a function from $\Gr$ to $\mathbb{C}.$
We denote the Fourier transform of $f$ by~$\FF{f},$
\begin{equation}\label{F:Fourier}
  \FF{f}(\xi) =  \sum_{x \in \Gr} f(x) e( -\xi \cdot x) \,,
\end{equation}
where $e(x) = e^{2\pi i x}$.
We rely on the following basic identities
\begin{equation}\label{F_Par}
    \sum_{x\in \Gr} f(x) \ov{g} (x)
        =
            \frac{1}{|\Gr|} \sum_{\xi \in \FF{\Gr}} \widehat{f} (\xi) \ov{\widehat{g} (\xi)} \,.
\end{equation}
\begin{equation}\label{svertka}
    \sum_{y\in \Gr} \Big|\sum_{x\in \Gr} f(x) g(y-x) \Big|^2
        = \frac{1}{|\Gr|} \sum_{\xi \in \FF{\Gr}} \big|\widehat{f} (\xi)\big|^2 \big|\widehat{g} (\xi)\big|^2 \,.
\end{equation}
If
$$
    (f*g) (x) := \sum_{y\in \Gr} f(y) g(x-y) \quad \mbox{ and } \quad
        (f\circ g) (x) := \sum_{y\in \Gr} f(y) g(y+x) = (g\c f) (-x)
$$
 then
\begin{equation}\label{f:F_svertka}
    \FF{f*g} = \FF{f} \FF{g} \quad \mbox{ and } \quad \FF{f \circ g} = \ov{\FF{\ov{f}}} \FF{g} = \FF{f}^c \FF{g} \,.
\end{equation}
Note also that
\begin{equation}\label{f:c*_and_*c}
    (f*g) (x) = (f^c \c g) (x) = (f \c g^c) (-x) \quad \mbox{ and } \quad
        (f \c g) (x) = (f^c * g) (x) = (f * g^c) (-x) \,,
\end{equation}
where for a function $f:\Gr \to \mathbb{C}$ we put $f^c (x):= f(-x)$.
Clearly,  $(f*g) (x) = (g*f) (x)$, $x\in \Gr$.
By $\langle f,g\rangle$ denote the scalar product of two complex functions $f$ and $g$.
Put $\langle f \rangle = \langle f,1 \rangle$, where $1$ is the constant function on $\Gr$.
We
will
write $\sum_x$ and $\sum_{\xi}$ instead of $\sum_{x\in \Gr}$ and $\sum_{\xi \in \FF{\Gr}}$ for simplicity.

We use in the paper  the same letter to denote a set
$S\subseteq \Gr$ and its characteristic function $S:\Gr\rightarrow \{0,1\}.$
By $|S|$ denote the cardinality of $S$.
Given $a\in \Gr$ we write $\d_a (x)$ for the delta function at the point $a$.
For a positive integer $n,$ we set $[n]=\{1,\ldots,n\}$.
All logarithms $\log$ are base $2$.
Signs $\ll$ and $\gg$ are the usual Vinogradov's symbols.

 For a sequence $s=(s_1,\dots, s_{k}) \in \Gr^k$ put $A_s=A\cap (A-s_1) \cap \dots \cap (A-s_{k}).$
Let
\begin{equation}\label{def:E_k}
    \E_k(A)=\sum_{x\in \Gr} (A\c A)(x)^k = \sum_{s_1,\dots,s_{k-1} \in \Gr} |A_s|^2 \,.
\end{equation}
If $k=2$ then $\E_{k} (A)$ is denoted by $\E(A)$
and is called the {\it additive energy} of $A$, see \cite{Tao_Vu_book}.
Some results on the quantities $\E_k (A)$ can be found in \cite{ss_E_k,sv}.
For any complex function $f$ and a positive integer $k$ denote by
$\Cf_{k+1} (f) (x_1,\dots,x_k)$ the quantity
$$
    \Cf_{k+1} (f) (x_1,\dots,x_k)
        = \sum_z f(z) f(z+x_1) \dots f (z+x_k) \,.
$$
The next lemma is a very special case of  Lemma 4 of paper \cite{s_ineq} and is the simplest generalization
of the second formula from (\ref{def:E_k}).

\Lemma
{\it
    Let $f,g$ be two complex functions on an abelian group $\Gr$.
    Suppose that $k$ is a positive integer.
    Then
$$
    \sum_{x_1,\dots,x_k} \Cf_{k+1} (f) (x_1,\dots,x_k) \Cf_{k+1} (g) (x_1,\dots,x_k)
        =
            \sum_z (f\c g)^{k+1} (z) \,.
$$
}
\label{l:E_k-identity}

\bigskip

Now consider the case $\Gr$ be the field.
If $q=p^s$, $p$ is a prime number then we write $\F_q$ for such finite field.
In the case $q=p$ by $R$ and $N$ denote the sets of quadratic residues and non--residues of $\F_p$,
correspondingly.
Clearly, $|R| = |N| = \frac{p-1}{2} := t$,
and $0\notin R$, $0\notin N$.
By $\chi_0$ denote the  principle character,
and denote the Legendre symbol on $\F_p$ by $\chi$.
Given nonzero $\la \in \F_p$ and a set $A \subseteq \F_p$ we will write
$$
    \la \cdot A := \{ \la \cdot a ~:~ a\in A \} \,.
$$
Thus $2A:=A+A \neq 2\cdot A$ in general.

\bigskip

\Def
{
Let $\_phi,\psi$ be two characters on $\F_q$.
The {\it Jacobi sum} $J(\chi,\psi)$ is defined by
$$
    J(\_phi,\psi) = \sum_x \_phi(x) \psi(1-x) \,.
$$
}

We need in a lemma (see \cite{BEW_book}, chapters 1--2).

\Lemma
{\it
    For any non principle character $\psi$, we have
    \begin{equation}\label{f:J_chi_chi_ov}
        J(\psi,\ov{\psi}) = - \psi(-1) \,.
    \end{equation}
    Whence
    \begin{equation}\label{f:K_chi_convolution}
        (\psi \c \ov{\psi}) (x) = p \d_0 (x) - 1 \,, \quad \quad  (\psi * \ov{\psi}) (x) = \chi(-1) (p \d_0 (x) - 1) \,.
    \end{equation}
    Further
\begin{displaymath}
G(p) := \sum_{x} \chi (x) e^{2\pi ix/p}
= \left\{ \begin{array}{ll}
\sqrt{p} &  \textrm{if ~ $p \equiv 1$ $mod$ $4$}\\
i\sqrt{p} & \textrm{if ~ $p \equiv -1$ $mod$ $4$}
\end{array} \right.
\end{displaymath}
    In particular, for any $x\neq 0$ the following holds $|\FF{R} (x)| \le \frac{\sqrt{p}+1}{2}$.
}
\label{l:K_chi}
\\
\Proof
Indeed, by the definition of Gauss sum, we have for all $x\in \F_p$ that
$$
    \FF{R} (x) = \frac{1}{2} \left( p \d_0 (x) - 1 + G(p) \chi(-x) \right)
$$
and the result follows.
$\hfill\Box$

Recall a well--known consequence  of the lemma.

\Lemma
{\it
    Let $g,h : \F_p \to \C$ be any complex functions.
    Then
\begin{equation}\label{f:rho_of_sum_1}
    \left| \sum_{x,y} g(x) h(y) \chi(x+y) \right|
        \le
            \| g \|_2 \left( p \| h \|^2_2 - |\langle h \rangle|^2 \right)^{1/2}
        \le
            \| g \|_2 \| h \|_2 \sqrt{p} \,,
\end{equation}
    and
\begin{equation}\label{f:rho_of_sum_2-}
    ((g\c \chi) \c (h\c \chi)) (x) = p (h\c g) (x) - \langle g \rangle \cdot \langle h \rangle \,.
\end{equation}
    In particular
\begin{equation}\label{f:rho_of_sum_2}
    \sum_{z} (g\c \chi) (x) \ov{(h\c \chi) (x)} = p \langle g,h \rangle - \langle g \rangle \cdot \langle \ov{h} \rangle \,.
\end{equation}
}
\label{l:rho_of_sum}

Note that inequality (\ref{f:rho_of_sum_1}) is sharp (see e.g. section \ref{sec:appendix}).
Formula (\ref{f:rho_of_sum_2}) of lemma above implies the "Cauchy--Schwartz"\, inequality in $\F_p$.

\Cor
{\it
    For any complex function $f : \F_p \to \C$, we have
\begin{equation}\label{}
    \| f \|_2^2 = \frac{|\langle f \rangle|^2}{p} + \frac{1}{p} \sum_x |(f\c \chi) (x)|^2 \ge \frac{|\langle f \rangle|^2}{p} \,.
\end{equation}
}
\label{cor:CS_new}

Also, using Lemma \ref{l:rho_of_sum} one can obtain simple upper bounds for the cardinalities
of sets $A,B$ such that $A+B\subseteq R$ or $A+B\subseteq N$
(see e.g. the proof of Theorem \ref{t:symmetric_difference} below or \cite{Bollobas,Palo_Alto_problems}).

Applying
Lemma \ref{l:E_k-identity}
we can easily improve Theorem \ref{t:Sarkozy} from \cite{Sarkozy_R}
(similar result was obtained in \cite{Shp_A+B}).

\Cor
{\it
    Let $A+B = R$.
    Then as $p\to \infty$, one has
    $$
        \left( \frac{1}{6} - o(1) \right)\sqrt{p} \le |A|,|B| \le (3+o(1)) \sqrt{p} \,.
    $$
}
\label{cor:Sarkozy}
\\
\Proof
We can assume that $A$ and $B$ are sufficiently large sets (see \cite{Sarkozy_R}).
By Lemma \ref{l:E_k-identity}, we have
$$
    |A|^4 |B| = \sum_{x\in B} (\chi \c A)^4 (x) \le \sum_{x} (\chi \c A)^4 (x)
        =
            \sum_{x,y,z} \Cf_4 (\chi) (x,y,z) \Cf_4 (A) (x,y,z) \,.
$$
Formula (\ref{f:Weil}) gives
$$
    |\Cf_4 (\chi) (x,y,z)| \le 3 \sqrt{p}
$$
with three exceptions : $x=y \neq 0, z=0$, further $x=z \neq 0, y=0$,
and $y=z \neq 0, x=0$.
Thus
$$
    |A|^4 |B| \le 3 \sqrt{p} |A|^4 + 3 p |A|^2
$$
and we obtain that $|A|, |B| \le \left( 3 + o(1) \right)\sqrt{p}$.
Because of $|A||B| \ge t = (p-1)/2$, we get $|A|,|B| \ge \left( \frac{1}{6} - o(1) \right)\sqrt{p}$.
This completes  the proof.
$\hfill\Box$

\bigskip

\bigskip


At the end of the section we recall the notion of perfect difference sets.

\Def
{
Let $\Gr$ be a group.
A set $A\subseteq \Gr$ is called a {\it perfect difference set}
if the convolution $(A\c A) (x)$ does not depend on the choice of $x\neq 0$.
}

Because of $(A\c A) (0) = |A|$
the definition above says that $(A\c A) (x) = (|A|-\la) \d_0 (x) + \la$,
where $\la \ge 1$ is some constant.
We will say that $A$ is
{\it $\la$--perfect difference set} in the case and just
{\it perfect difference set}  for $\la = 1$.
Generally speaking, let $\mathcal{L}$ be the algebra of the functions of the form $a \d_0 (x) + b$,
where $a$ and $b$ are some complex  constants.
Clearly, a function $\_phi$ belongs to $\mathcal{L}$
iff $\FF{\_phi} \in \mathcal{L}$.
Denote by $\mathcal{D}$ the set of all {\it perfect difference functions},
that is all functions $\_phi$ such that $\_phi \c \_phi$ or, equivalently, $\FF{\_phi}^c \cdot \FF{\_phi}$
belongs to $\mathcal{L}$.
For example, Lemma \ref{l:K_chi} says that $\chi \in \mathcal{D}$.
It is easy to check that $\mathcal{D}$ is closed under convolutions $*$ and $\c$.
Moreover, if $\_phi \in \mathcal{D}$, $\FF{\_phi} (x) \neq 0$, $x\in \Gr$
and $\_phi \c \psi$ or $\_phi * \psi \in \mathcal{D}$
then $\psi \in \mathcal{D}$.
Finally, if $f\in \mathcal{D}$ then $\a+\beta f \in \mathcal{D}$ for any $\a,\beta \in \C$.

\bigskip

We need in a simple lemma.

\Lemma
{\it
    For all $x\neq 0$, we have
    \begin{equation}\label{f:R_conv}
        (R\circ R) (x) = \frac{p-3}{4} - \frac{\chi(x)}{4} \left(1 + \chi(-1) \right) \,.
    \end{equation}
}
\label{l:R_conv}
\Proof
Since
$$
    R(x) = \frac{1}{2} \left( \chi_0 (x) + \chi(x) \right) 
$$
it follows that for $x\neq 0$
$$
    (R\circ R) (x) = \frac{1}{4} \sum_z (\chi_0 (z) + \chi(z)) (\chi_0 (z+x) + \chi(z+x))
        =
$$
$$
        =
            \frac{p-2}{4} - \frac{1}{4} \left( \chi(x) + \chi(-x) \right) + \frac{1}{4} \sum_z \chi(z(z+x))
$$
and the result follows from Lemma \ref{l:K_chi}.
$\hfill\Box$

\bigskip

In particular, if $p \equiv -1 \pmod 4$ then $R\in \mathcal{D}$.

\bigskip

Let $A \subseteq \Z_P$ be a perfect difference set.
A residue  $m$ is called a {\it multiplier} of $A$ if $mA=A$.
We formulate the multiplier theorem, see e.g. \cite{Hall}.

\Th
{\it
    Let $A$ be a $\la$--perfect difference set, $A\subseteq \Z_P$.
    Suppose that $m$ be any prime number such that $m | (|A|-\la)$, $(m,P)=1$
    and $m>\la$.
    Then $m$ is a multiplier of some translation of $A$.
}
\label{t:multiplier}

Clearly, the set of multipliers of $A$ forms a group and,
moreover,
one can choose a translation of $A$ fixing by the group (see \cite{Hall}).

Recall also a beautiful
theorem of Singer (see e.g. \cite{Hall})
concerning finite projecting geometries.


\Th
{\it
    Suppose that $P$ is a number of the form $P=n^2+n+1$, $n=q^s$, where $s\ge 1$ and $q$ is a prime number.
    Then there is a perfect difference set $A\subseteq \Z_P$ such that $(A\c A) (x) =1$ for all $x\neq 0$.
}
\label{t:Singer}

\refstepcounter{section}
\label{sec:proof}

\begin{center}
{\bf \large
    \arabic{section}. Sumsets and differences}
\end{center}

We begin with a simple lemma.

\Lemma
{\it
    Let $c$ be
    an
    integer, and $\psi : \Gr \to \Z$ be a function.
    Then
    \begin{equation}\label{f:div_c}
        \| \psi \|_2^2 \ge c|\sum_x \psi (x) | - (c-1) \cdot|\sum_{x ~:~ 0 < |\psi (x)| < c } \psi (x)| \,.
    \end{equation}
    Further
    \begin{equation}\label{f:div_c'}
        \| \psi \|_2^2 = c\sum_x \psi (x) + \sum_{k} |\{ x ~:~ \psi(x) = k \}| \cdot (k^2-ck) \,.
    \end{equation}
}
\label{l:div_c}
\Proof
Let $\sigma = \sum_x \psi (x)$.
We can suppose that $\sigma >0$, otherwise consider the function $-\psi$.
For any integer $k$ put
$$
    p_k = |\{ x ~:~ \psi(x) = k \} | \,.
$$
Then
$$
    \sigma = \sum_{k} k p_{k} = \sum_{k ~:~ |k| \ge c} k p_{k} + \sum_{k ~:~ 0 < |k| < c} k p_{k}  = \sigma_1 + \sigma_2 \,,
$$
and hence
$$
    \| \psi \|_2^2 = \sum_{k} k^2 p_{k}
        =
            \sum_{k ~:~ |k| \ge c} k^2 p_{k} + \sum_{k ~:~ 0 < |k| < c} k^2 p_{k}
                \ge
        c\sigma_1 + \sigma_2
            =
                c\sigma - (c-1) \sigma_2
                    \ge
$$
$$
                    \ge
                        c\sigma - (c-1) |\sigma_2|
$$
as required.
Formula
(\ref{f:div_c'})
follows similarly.
This completes the proof.
$\hfill\Box$

\bigskip

Now we can prove a small generalization of the first part of our  main Theorem \ref{t:main_intr}.

\bigskip

\Th
{\it
    Let $p$ be a prime number, $R \subseteq \F_p$ be the set of quadratic residues and $A\subseteq \F_p$ be a set.
    If $A+A = R$ then $p=3$ and  $A=\{ 2 \}$.
    Moreover for all sufficiently large $p$, we have
    \begin{equation}\label{f:symmetric_difference}
        \max\left\{ |R\setminus (A+A)|, \sum_{x\in (A+A) \setminus R} (A*A) (x) \right\}
            \ge
                \left( \frac{1}{6} - o(1) \right) |A| \,.
    \end{equation}
}
\label{t:symmetric_difference}
\Proof
Suppose that $A,B$ are two sets such that $A+B\subseteq R$.
Let $a=|A|$, $b=|B|$.
Define the function $\eps(x)$ by the formula
\begin{equation}\label{f:eps_def_I}
    ( A \c \chi ) (x) = a B(x) + \eps (x) \,.
\end{equation}
We have $\eps (x) = 0$, $x\in B$.
Using
Lemma \ref{l:rho_of_sum}, we get
\begin{equation}\label{f:eps_l_2_pred}
    \| \eps \|_2^2 = pa-a^2 - a^2 b \,.
\end{equation}
Further
\begin{equation}\label{f:eps_average}
    \langle \eps \rangle = -ab \,,
\end{equation}
thus, by the Cauchy--Schwartz inequality
$$
    a^2 b^2 \le (p-b)(pa-a^2 - a^2 b)
$$
or
\begin{equation}\label{f:1_estimate}
    p + \frac{ab}{p} \ge ab+a+b \,.
\end{equation}
On the other hand, by formula (\ref{f:rho_of_sum_1}) of Lemma \ref{l:rho_of_sum}
or just by inequality (\ref{f:eps_l_2_pred}),
one has $ab<p$.

Now suppose that $A=B$ and $A+A=R$.
If $a=1$ then $\frac{p-1}{2}=1$, whence, $p=3$ and, clearly, $A=\{ 2 \}$.
Thus, suppose that $a>1$.
Because of $A+A=R$, we
obtain
$\binom{a}{2}+a \ge \frac{p-1}{2}$ or, in other words, $a^2 + a \ge p-1$.
Using the last estimate, the fact $a^2 < p$ and  inequality (\ref{f:1_estimate}), we get
$$
    p+1 > p + \frac{a^2}{p} \ge a^2+2a \ge p-1+a
$$
with contradiction.

Now, let us prove  the second part of the theorem.
Put
$$
    |R\setminus (A+A)| := \zeta_1 a \quad \mbox{ and } \quad
        \sum_{x\in (A+A) \setminus R} (A*A) (x) := \zeta_2 a \,.
$$
Let also $\zeta = \max \{ \zeta_1, \zeta_2 \}$.
We need to obtain the lower bound for $\zeta$.
Suppose that $\zeta \le \frac{1}{6} + o(1)$, $p\to \infty$.
We have
$$
    \binom{a}{2}+a \ge \frac{p-1}{2} - \zeta_1 a \ge \frac{p-1}{2} - \zeta a
$$
or, in other words,
\begin{equation}\label{tmp:09.06.2012_1}
    a^2 + a(1+2\zeta) \ge p-1 \,.
\end{equation}
Put
$$
    (A\c \chi) (x) = a A(x) + \eps (x) \,.
$$
As in (\ref{f:eps_average}) the average value of $\eps(x)$ equals $-a^2$,
and we estimate $l_2$--norm of the function as in (\ref{f:eps_l_2_pred})
$$
    \| \eps \|_2^2 = \sum_{x} ( (A\c \chi) (x) - a A(x))^2
        =
            pa-a^2 + a^3 - 2a \sum_{x\in A} (A\c \chi) (x)
                =
$$
$$
                =
                     pa-a^2 - a^3 + 2a \sum_x (1-\chi (x)) (A*A) (x)
                        =
$$
\begin{equation}\label{f:eps_l_2}
                        =
                            pa-a^2 - a^3 + 2a \sum_{x\in (A+A) \setminus R} (1-\chi (x)) (A*A) (x)
                                \le
                                    pa-a^2 - a^3 + 4 \zeta a^2 \,.
\end{equation}
By the Cauchy--Schwartz inequality
\begin{equation}\label{tmp:15.05.2013_1}
    \frac{a^4}{p} \le \| \eps \|_2^2 \le pa-a^2 - a^3 + 4 \zeta a^2 \,.
\end{equation}
Note that from the last estimate it follows that $a\ll \sqrt{p}$.
Using (\ref{tmp:15.05.2013_1}) and inequality (\ref{tmp:09.06.2012_1}),
we obtain after some calculations that $\zeta \ge \frac{1}{6} - o(1)$.
This completes the proof.
$\hfill\Box$

\bigskip

Now we
reprove
a result of C. Bachoc, M. Matolcsi, I.Z. Rusza from \cite{Ruzsa_sqrt}.

\bigskip

\Th
{\it
    Let $p$ be a prime number, $R \subseteq \F_p$ be the set of quadratic residues and $A\subseteq \F_p$ be a set.
    If $A-A \subseteq  R \sqcup \{ 0 \}$ then
\begin{displaymath}
p \ge
    \left\{ \begin{array}{ll}
        |A|^2 + |A| - 1 & \textrm{if ~ $|A|$ ~ is even} \\
        |A|^2 + 2|A| - 2& \textrm{if ~ $|A|$ ~ is odd.}
    \end{array} \right.
\end{displaymath}
}
\label{t:Lev_q}
\\
\Proof
Put $|A| = a$.
Clearly, $p\equiv 1 \pmod 4$.
As in the previous theorem, we have
\begin{equation}\label{tmp:10.05.2013_0}
    ( A * \chi ) (x) = (a-1) A(x) + \eps (x) \,,
\end{equation}
where $\eps (x) = 0$ for $x\in A$.
As above
$$
    \| \eps \|_2^2 = pa-a^2-a(a-1)^2 \,, \quad \mbox{ and } \quad \langle \eps \rangle = -a(a-1) \,.
$$
Further
\begin{equation}\label{f:eps_eo}
    \eps(x) = |R\cap (x-A)| - |N\cap (x-A)| = 2 |R\cap (x-A)| - a
\end{equation}
for all $x\notin A$.
First of all consider the case $a$ is an even number.
By formula (\ref{f:eps_eo}) the values of $\eps(x)$ are even for all $x$
(recall that $\eps(x) = 0$ on $A$).
Using Lemma \ref{l:div_c} with $c=2$, we obtain
$$
    \| \eps \|_2^2 = pa-a^2-a(a-1)^2  \ge 2a(a-1)
$$
or, in other words,
$$
    p \ge a^2+a-1
$$
as required.

Now suppose that $a$ is odd and $a>1$.
Of course $\eps(x) = 0$ on $A$, and by formula (\ref{f:eps_eo}) all other values
of $\eps(x)$ are odd.
First of all note that there is $x\notin A$ such that $\eps (x) \neq -1$.
Indeed, otherwise, applying (\ref{tmp:10.05.2013_0}), we get
\begin{equation}\label{tmp:10.05.2013_1}
    \langle A, \chi \rangle = \langle \eps, \chi \rangle = p A(0) - a - (a-1) \langle A,\chi \rangle \,.
\end{equation}
Shifting, we can suppose that $0\in A$ and, hence,
$\langle A,\chi \rangle = a-1$.
Substituting into (\ref{tmp:10.05.2013_1}) gives $p=a^2$ with contradiction.
(Note that nevertheless it can be the case for the field $\F_q$ with $q$ equals a prime power,
see \cite{Ruzsa_sqrt}).
Now give an another proof of the fact.
Put $p_k := |\{ x \in \F_p ~:~ \eps(x) = k \}|$, $k\in \Z$.
Clearly, $p_0 = a$ and $p_{2l} = 0$ for all nonzero integer $l$.
As before, using formula (\ref{f:div_c'}) of Lemma \ref{f:div_c} with $c=-2$
and the fact that
the quantity $p_{-1}$
does not exceed $(p-1)/2$
(consider $\eps(x) \pmod p$ and note that $\chi (x) \equiv x^{(p-1)/2} \pmod p$), we obtain
$$
    \| \eps \|_2^2 = pa-a^2-a(a-1)^2 = -2 \langle \eps \rangle + \sum_{k} (k^2+2k) p_k
        \ge
            2 a(a-1) - p_{-1} + 3 \sum_{k \neq -1,0} p_k
                =
$$
$$
                =
    2 a(a-1) + 3p - 4p_{-1} - 3a
        \ge
            2 a(a-1) + p + 2 - 3a
$$
or, in other words,
$$
    p \ge a^2+a-4 + \frac{p+2}{a} \,.
$$
After some calculations we
get $p \ge a^2 + 2a-2$, provided by $a>1$.
This completes the proof.
$\hfill\Box$

\refstepcounter{section}
\label{sec:proof2}

\begin{center}
{\bf \large
    \arabic{section}. Restricted sumsets}
\end{center}

The aim of the section is to prove the second part of our main Theorem \ref{t:main_intr}.

Write any prime
$p$
as $p=n^2 + \Delta$, $1\le \Delta \le 2n$.
We are interested in sets $A\subseteq \F_p$ such that $A\hat{+} A = R$.
Let us consider some examples.

\bigskip

\Exm
{
There are
examples of {\it perfect difference sets} $A \subseteq \F_p$ with $A\hat{+} A = R$ : \\
{\bf (a)}$~$ $p=3$, $A=\{ 0,1 \}$.\\
{\bf (b)}$~$ $p=7$, $A=\{ 3,5,6 \}$.\\
{\bf (c)}$~$ $p=13$, $A=\{ 0,1,3,9 \}$, $A=\{ 0,4,10,12 \}$.\\
Note that in all these examples $p=n^2+n+1$, $|A|=n+1$,
$n=q^s$, $q$ is a prime number or $1$.
In particular $\Delta = n+1$.
The existence of a perfect difference set in $\F_p$ for $p$ of such a  form
guaranties   by Singer's theorem \ref{t:Singer}.
Note also, that it is easy to check
using the multiplier theorem and another tools
(see the beginning of the  proof of Proposition \ref{pr:pd_sets}
and Theorem \ref{t:A+A=R}
below)
that the sets above form the complete list of
(perfect difference) sets
$A \subseteq \F_p$ with $A\hat{+} A = R$ for
$p=3,7,13$.
}
\label{exm:pd_sets}

\bigskip

In our proof of the second part of main Theorem \ref{t:main_intr}
we begin with the situation when $A$ is a perfect difference sets.
Surprisingly, it turns out that it is the most important  case.

\bigskip

\Pred
{\it
    Let $p$ be a prime number, $p=n^2+n+1$,
    $R \subseteq \F_p$ be the set of quadratic residues
    and $A\subseteq \F_p$
    be a perfect difference set.
    Suppose that $A\hat{+} A = R$.
    Then the set $A$ must be
    a set
    from Example \ref{exm:pd_sets}.
}
\label{pr:pd_sets}
\\
\Proof
Let $a=|A|$.
Because of $A$ is a perfect difference set, we have $a=n+1$.
Firstly,
$(A\hat{+} A) \cap 2\cdot A = \emptyset$
and, consequently, $R\cap 2\cdot A = \emptyset$.
In particular either $A\subseteq R\sqcup \{0\}$ or $A\subseteq N$.
Secondly,
by Theorem \ref{t:multiplier} there is $s \in \F_p$ such that the set $A+s$ is fixed by a multiplier $m \neq 1$.
Then $(A+s)\hat{+} (A+s) = R+2s$ and hence $m(R+2s) = mR +2sm = R+2s$.
If $m\in R$ we have automatically $s=0$.
If $m\in N$ then $mR=N$, and hence
$N(x) = R(x+2s(m-1))$.
In other words
$$
    \chi_0(x) - \chi(x) = \chi_0(x+2s(m-1)) + \chi(x+2s(m-1)) \,.
$$
Clearly, the case $s= 0$ is impossible.
Multiplying by $\chi(x)$, summing and using Lemma \ref{l:K_chi}, we get
$$
    p-1=1 - \chi(-2s(m-1)) \,.
$$
That implies $p=1$ or $p=3$.
The case $p=3$ is Example \ref{exm:pd_sets} (a).
Thus for $p>3$, we obtain that any multiplier of $A$ belongs to $R$.
Below we assume that $n>1$ and, hence, $p>3$.
By Theorem \ref{t:multiplier} and the previous arguments any prime factor of $n$ is a quadratic residue,
which implies that $n$ belongs to $R$ itself.
We have $p=n^2+n+1$, $n^3 \equiv 1 \pmod p$ and hence the order of $n$ is three
(in particular $p\equiv 1 \pmod 3$).
The same follows from the fact that
$n\equiv 0 \pmod 3$ if $0\in A$
and
$n\equiv -1 \pmod 3$ if $0\notin A$.
Take an arbitrary nonzero $x \in A$.
Because of $n>1$, we see that  $x,xn,xn^2$ are different and belong to $A$.
Whence $x+xn^2 \equiv x -(n+1)x \equiv -nx \in R$, so $x\in -R$.
That is equivalent $A \setminus \{ 0 \} \subseteq -R$.

Our arguments rest on the identity
\begin{equation}\label{tmp:05.05.2012_1}
    R(x)
    = \frac{1}{2} ( (A*A) (x) - (2\cdot A) (x) )
\end{equation}
which is a consequence of the fact that
$A$ is a perfect difference set.

First of all suppose that
$A\subseteq N$.
Then $A\subseteq -R$, hence, $-1\in N$, $p\equiv -1 \pmod 4$ and $2\in R$.
Since $R\in \mathcal{D}$
and, simultaneously, $A\in \mathcal{D}$
it follows from
formula (\ref{tmp:05.05.2012_1})
that
$((A*A) \c 2\cdot A) (x) + (2\cdot A \c (A*A)) (x) \in \mathcal{L}$.
Using identity (\ref{tmp:05.05.2012_1}) again, we see that
$(R\c (2\cdot A)) (x) + ((2\cdot A) \c R) (x) \in \mathcal{L}$.
It follows that
\begin{equation}\label{tmp:05.05.2012_2}
    | R\cap (2\cdot A + x)| + |R\cap (2\cdot A - x)| = a \,, \quad \forall x \neq 0 \,.
\end{equation}
We know the number
$2$ belongs to $R$.
If $x\in 2\cdot A$ then
$$
    | R\cap (2\cdot A + x)| = |R\cap (A+x/2)| = a-1
$$
because of $x/2\in A$ and $A\hat{+} A = R$.
On the other hand,
using the fact
$2\in R$ again, we
get
\begin{equation}\label{tmp:05.05.2012_3}
    \sum_{x\in 2\cdot A} |R\cap (2\cdot A - x)| = \sum_x R(x) (A\c A) (x) = t
\end{equation}
and hence there is $x\in 2\cdot A$ such that $|R\cap (2\cdot A - x)| \ge t/a$.
Combining the last inequality with (\ref{tmp:05.05.2012_2}), we obtain
$$
    a-1 + \frac{t}{a} \le a
$$
or $a\le 3$.
That is Example \ref{exm:pd_sets} (b).

It remains consider the situation $A\subseteq R \sqcup \{ 0 \}$.
In the case $p\equiv 1 \pmod 4$ and $2\in N$ but, unfortunately,  $R\notin \mathcal{D}$,
and, thus, we need in more delicate arguments.
Using (\ref{tmp:05.05.2012_1}),
one can see that
$$
    (R\c R) (x) = \frac{1}{4} \big( ((A*A) \c (A*A))(x) + ((2\cdot A) \c (2\cdot A)) (x)
        -  ((2R+2\cdot A) \c 2\cdot A) (x)
            -
$$
$$
            -
                2\cdot A \c (2R+2\cdot A) \big)
                    =
                        \frac{1}{4} \big( ((A*A) \c (A*A))(x) - ((2\cdot A) \c (2\cdot A)) (x)
                            -
                                2 (R\c 2\cdot A) (x) - 2 (2\cdot A \c R) (x) \big) \,.
$$
Because of $\E(A) = 2a^2 -a$, and $A*A \in \mathcal{D}$, we have
$$
    ((A*A) \c (A*A))(x) = \frac{a^4-\E(A)}{p-1} = \frac{a^4-2a^2+a}{p-1} = a^2+a-1 \,, \quad \forall x \neq 0 \,.
$$
Combining the last formula, Lemma \ref{l:R_conv}, and the fact $A\in \mathcal{D}$, we obtain
$$
    p-3-2\chi(x) = a^2+a-1 - 1 - 2 (R\c 2\cdot A) (x) - 2 (2\cdot A \c R) (x) \,, \quad \forall x \neq 0 \,.
$$
In other words
\begin{equation}\label{tmp:05.05.2012_4}
    a+\chi(x) = (R\c 2\cdot A) (x) + (2\cdot A \c R) (x) \,, \quad \forall x \neq 0 \,.
\end{equation}
Let $A^* = A\setminus \{ 0 \}$.
Using the arguments similar (\ref{tmp:05.05.2012_2})---(\ref{tmp:05.05.2012_3}),
we see that for all $x\in 2\cdot A^*$ the following holds
$| R\cap (2\cdot A + x)| = |N\cap (A+x/2)| = 1$
and
$$
    \sum_{x\in 2\cdot A^*} |R\cap (2\cdot A - x)| = \sum_x N(x) (A^*\c A) (x) =
        \sum_x N(x) (A\c A) (x) - |A\cap N| = t \,.
$$
Whence there is $x\in 2\cdot A^*$ such that $|R\cap (2\cdot A - x)| \le t/(a-1)$.
For any $x\in 2\cdot A^*$, we have $\chi(x)=-1$.
Applying (\ref{tmp:05.05.2012_4}), we obtain
$$
    a-1 \le 1+ \frac{t}{a-1}
$$
or $a\le 4$ (if $0\notin A$ then it is easy to check similarly  that $a\le 3$).
That is Example \ref{exm:pd_sets} (c).
Indeed, $A\subseteq R\sqcup \{0\}$,
$R= \{ 1,3,4,9,10,12 \}$ and it is easy to
see
that $A$ either
$\{ 0,1,3,9 \}$ or $\{ 0,4,10,12 \}$
(note that $|A|-1=3$ is a multiplier of $A$).
%
%
%
This completes the proof.
$\hfill\Box$

\bigskip

Note that Example \ref{exm:pd_sets} shows the possibility of the
cases $\Delta=2,3$ as well as $\Delta=n+1$.

\bigskip

Now we can consider the general situation.
Of course, the third part of theorem below is the main one but it follows
from the first.
The second part shows that $A$ with $A\hat{+} A = R$ is close to a perfect difference set in some sense.

\bigskip

\Th
{\it
    Let $p$ be a prime number, $p=n^2 + \D$, $1\le \D \le 2n$,
    $R \subseteq \F_p$ be the set of quadratic residues and $A\subseteq \F_p$ be a set.
    \\
    ${\bf 1)~}$ If $A\hat{+} A = R$ and $A$ is not from Example \ref{exm:pd_sets} (a)
    then $|A|=n+1$, $3\le \Delta \le n+1$ and
    $|2\cdot A \cap R| \le \frac{\sqrt{|A| \Delta - 3|A|+1}}{2}$.
    If $0\notin A$ then $\chi(2)=1$, and $|A| \le 6$.
    If $0\in A$ then $\chi(2)=-1$ and if $\D = n+1$ then $A$ is a perfect difference set.
    \\
    ${\bf 2)~}$ If $A\hat{+} A = R$ then it is close to a perfect difference set
    in the sense
    \begin{equation}\label{f:A_energy_close}
        \E (A) = 2|A|^2 - |A| + \mathcal{E}_1 \,,
    \end{equation}
    where
    $$
        \mathcal{E}_1 \le 6(|A|-\Delta) + 2|A| + \min\{ (|A|-\Delta) |A|^{}, |A| \sqrt{3 (|A|-\D)}, (|A|-\Delta)^2 \} \,,
    $$
    and
    \begin{equation}\label{f:A_difference_close}
        (A\c A) (x) = (|A|-1) \d_0 (x) + 1 + \mathcal{E}_2 (x) \,,
    \end{equation}
    where $\sum_x \mathcal{E}_2 (x) = |A| - \Delta$, $\| \mathcal{E}_2 \|^2_2 = \mathcal{E}_1 + \Delta - |A|$.
    \\
    ${\bf 3)~}$ If $A\hat{+} A = R$ then $A$ is a perfect difference set such that
    if $|A|$ is even then $|A| \le 6$ and if $|A|$ is odd then $|A| \le 5$.
}
\label{t:A+A=R}
\\
\Proof ${\bf 1)~}$
Put $a=|A|$.
The assumption $A\hat{+} A = R$ implies that $\binom{a}{2} \ge \frac{p-1}{2}$.
In other words, $a^2 - a \ge p-1$.
The last inequality
is equivalent to
$a\ge n+1$.
We can assume that $a\ge 3$ because otherwise we have Example \ref{exm:pd_sets} (a).
One can also check that for $a\ge 3$ the case $p=5$ is impossible,
and hence we will assume that $p\ge 7$.
Below we will use the fact that $A\hat{+} A \subseteq R$ only.
Put
$$
    d := \sum_{x\in A} (\chi(2x) - 1) = \eta a \,.
$$
Clearly, $\eta \in [-2,0]$.
Let us obtain some further estimates on $d$ and $\eta$.
By Lemma \ref{l:rho_of_sum}, we get
\begin{equation}\label{tmp:30.05.2012_1}
    \sum_x (A \c \chi)^2 (x) = pa-a^2 \,.
\end{equation}
Further, similar to (\ref{f:eps_def_I}),
we have
\begin{equation}\label{f:decomposition_2A*}
    (A\c \chi) (x) = (a-1+\chi(2x)) A(x) + \eps (x)
\end{equation}
where $\eps(x) = 0$ for all $x\in A$.
Hence
$$
    a(a-1)^2 + 2(a-1)(d+a) + a-\o = \sum_{x\in A} (a-1 + \chi(2x))^2 \le pa-a^2 \,,
$$
where $\o = 1$ if $0\in A$ and $\o=0$ otherwise.
It is easy to calculate that the last inequality
and the bound $a^2 - a \ge p-1$ imply that
\begin{equation}\label{f:d_zero}
     d \le -a - \frac{a-\o}{2(a-1)} < -a \,.
\end{equation}
In particular, $\eta \in [-2,-1)$.
In the case $0\notin A$,
we have
$$
    (d+a)^2 + a(a-1)^2 + 2(a-1)(d+a) + a  \le pa-a^2
$$
and after some manipulations, we obtain
\begin{equation}\label{f:d_nonzero-}
    d\le - 2a+1 + \sqrt{a^2 - 3a+1} \,,
\end{equation}
and
\begin{equation}\label{f:d_nonzero}
    d\le - 2a+1 + \sqrt{a\Delta - 3a+1} \,,
\end{equation}
provided by $p=(a-1)^2 + \D$.
Note that in the case $\D \ge 3$.
Now suppose that $0\in A$. Write $A=A^* \bigsqcup \{ 0 \}$.
Then $A^* \subseteq R$, and the set $2\cdot A^*$ belongs either $R$ or $N$.
In the first case $d=-1$ and that is impossible by (\ref{f:d_zero}).
If $2\cdot A^* \subseteq N$ then $d=-2a+1$, $\chi (2) = -1$, and (\ref{f:d_nonzero-}), (\ref{f:d_nonzero}) take place.
Thus in any case
\begin{equation}\label{f:2A&R-}
    |2\cdot A \cap R| \le \frac{\sqrt{a^2 - 3a+1}}{2} \,,
\end{equation}
and
\begin{equation}\label{f:2A&R}
    |2\cdot A \cap R| \le \frac{\sqrt{a\Delta - 3a+1}}{2} \,,
\end{equation}
provided by $p=(a-1)^2 + \D$.
One of the main ideas of
the
further proof is an exploitation of inequalities (\ref{f:d_zero}), (\ref{f:2A&R-}), (\ref{f:2A&R}).
They mean that the distribution of the intersections $|R\cap (A+x)|$, $x\in \F_p$ is an asymmetric a little bit.

Returning to (\ref{f:decomposition_2A*}), we get
for all $x\notin A$ that
\begin{equation}\label{f:eps_not_A}
    \eps (x) = |R\cap (A+x)| - |N\cap (A+x)| = 2|R\cap (A+x)| - a + A(-x) \,.
\end{equation}
Clearly,
\begin{equation}\label{f:eps_av_dot}
    \langle \eps \rangle = - a^2 - d \,.
\end{equation}
As above, we have
$$
    \sum_{x\in A} (a-1+\chi(2x))^2 = a(a-1)^2 + 2(a-1)(d+a) + a - \o = a^3 + 2(a-1)d - \o \,,
$$
where, as before, the quantity $\o$ equals $|A\cap \{ 0 \}|$.
As in Theorem \ref{t:symmetric_difference}, we obtain an analog of (\ref{f:eps_l_2_pred})
\begin{equation}\label{f:eps_L_2_new}
    \| \eps \|_2^2 = pa-a^2 -a^3 - 2(a-1)d + \o \,.
\end{equation}
Thus by the Cauchy--Schwartz inequality and identity (\ref{f:eps_av_dot}), we get
\begin{equation}\label{tmp:08.03.2013_1}
    \frac{(a^2+d)^2}{p-a} \le pa-a^2 -a^3 - 2(a-1)d + \o \,.
\end{equation}
In other words
\begin{equation}\label{f:diskr}
   a^2 (p-1) + a b(\eta) - c(\eta) := a^2 (p-1) + a (\eta^2 + 2\eta + 2\eta p + 2p) - 2\eta p - p^2 - \o \frac{p-a}{a} \le 0 \,.
\end{equation}
It is easy to check that the formula for the right root of
(\ref{f:diskr}), namely,
$$
    x(\eta) = \frac{-b(\eta) + \sqrt{b^2(\eta) + 4(p-1) c(\eta)}}{2(p-1)}
$$
contains an decreasing function $-b(\eta)$ and an increasing function
$$
    g(\eta) := b^2(\eta) + 4(p-1) c(\eta)
    =
$$
$$
    = \eta^4+(4+4p)\eta^3+(4+4p^2+12p)\eta^2+16 \eta p^2+4p^3 + \o \frac{4(p-1)(p-a)}{a} \,.
$$
on $\eta \in [-2,-1)$.
Indeed, one can see that $g'(-2)=0$ and $g''$ grows on $\eta \in [-2,-1)$.
Put $e= \o \frac{4(p-1)(p-a)}{a}$.
Because of $a\ge \sqrt{p}$, we have $e\le 4p^{3/2}$.
Hence, we should substitute $\eta=-2$ into $b(\eta)$,
$\eta=-1$ into $g(\eta)$, and check the following inequality
\begin{equation}\label{tmp:12.05.2013_1}
    a\le \frac{2p + \sqrt{4p^3-12p^2+8p+1+e}}{2(p-1)} \le \sqrt{p}+1 < n+2 \,.
\end{equation}
One can insure  that, indeed, the second bound from (\ref{tmp:12.05.2013_1}) holds,
provided by $p\ge 7$.
Thus, for any $\eta$, we have $a=n+1$.
Recalling $a^2 - a \ge p-1$, we obtain
\begin{equation}\label{f:C_n^2_fake}
    a^2 - a = n^2 + n \ge p-1 = n^2 + \Delta -1 \,.
\end{equation}
Whence $3 \le \Delta \le n+1$ and in the case $\Delta = n+1$
all sums $a'+a''$ are different for different $a',a'' \in A$.
Using, the arguments above, we see that if $0\in A$ and $\D = n+1$ then $A$ is a perfect difference set.

Finally,
consider the case $0\notin A$.
We will prove, in particular, that $\chi(2)=1$ in the situation.
By (\ref{f:decomposition_2A*})
one has
$$
    \eps (0) = \sum_{x\in A} \chi(x) = \chi(2) \sum_{x\in A} \chi(2x) = \chi(2) (d+a) \,.
$$
Further,
identity
(\ref{f:decomposition_2A*}) and the last formula imply
$$
    \sum_x \eps(x) \chi(x) = -a - \chi(2) (a^2+da-d) \,.
$$
Hence, as in (\ref{tmp:08.03.2013_1}), excluding zero point from the function $\eps(x)$,
and, additionally, using Corollary \ref{cor:CS_new}, we get
\begin{equation}\label{tmp:09.03.2013_1}
        \frac{(a^2+d + \chi(2) (d+a))^2 + (a+\chi(2) (a^2+da-d))^2}{p}
    \le pa-a^2 -a^3 - 2(a-1)d - (d+a)^2 \,.
\end{equation}
After some manipulations, we obtain for $\chi(2)=-1$ that
$$
    h_{a,d} (\D) := a\D^2+(a^3-6a^2+2a-4da+2d-d^2)\D-
$$
\begin{equation}\label{tmp:10.03.2013_1'}
    -6a^4+13a^3-8a^2+14a^2 d-2d^2 a^2-6da^3+a-10da+4d^2 a+2d-2d^2
        \ge 0
\end{equation}
the maximum of $h_{a,d} (\D)$ is attained at $\D=a$ and is equal to
$$
    -(2a^2-3a+2)d^2 - (6a^3-10a^2+8a-2)d-5a^4+8a^3-6a^2+a \,.
$$
The integer maximum of the last expression is attained at $d \in [-1.5 a,-1.5a+1]$.
Substituting the integer values of $d$ from the interval, we obtain a contradiction.
Almost the same is true in the case $\chi(2)=1$.
In the situation
$$
    h^*_{a,d} (\D) := a\D^2+(a^3-6a^2+2a-4da+2d-d^2)\D
    +
$$
\begin{equation}\label{tmp:10.03.2013_1}
    +
        a-2d^2 a^2-6d^2+5a^3-8a^2+2d-6a^4+4d^2 a-10da+6a^2 d-6da^3 \ge 0
\end{equation}
Again the maximum of $h^*_{a,d} (\D)$
is attained at $\D=a$ and is equal to
$$
    -(2a^2-3a+6)d^2 - (6a^3 -2a^2 +8a-2)d-5a^4-6a^2+a \ge 0 \,.
$$
Here the integer maximum
is attained at $d\in (-1.5 a-3, -1.5a]$.
Substituting the integer values of $d$ from the interval into the last inequality, we obtain a contradiction for $a\ge 7$.
Thus, we have proved part ${1)}$ completely.
Let us make some additional remarks which we will use later.
If $a=3$ then $p=7$, $\chi(2) = 1$, $\D=a$, $|A\hat{+} A| = \binom{3}{2}=3$,
$R=\{1,2,4\}$.
If $A\cap R = \emptyset$ then $A$ is from Example \ref{exm:pd_sets} (b).
Further, it is easy to see that $0\notin A$, and $A$ can intersect $R$ just at one point.
If $1\in A$ then we have a contradiction because $5,6$ cannot belong to $A$ in the case.
The same is true for $2\in A$ and $4\in A$.
Thus $A$ is from Example \ref{exm:pd_sets} (b), provided by $a=3$.
Finally, if $a=5$ then $p=19$ because of $\D \ge 3$.
Hence $\chi(2)=-1$, $0\in A$, $d=-2a+1 = -9$.
Substituting this into (\ref{tmp:08.03.2013_1}), we obtain a contradiction.
The remain cases $a=4,6$ as well as the situation $0\in A$ will be considered later
in the third part of the proof.

${\bf 2)~}$
The fact $R=A\hat{+} A$ implies that
\begin{equation}\label{tmp:22.05.2012_1}
    (A*A) (x) = 2R(x) + (2\cdot A) (x) + Z (x) \,,
\end{equation}
where $Z(x) \ge 0$.
It is easy to see  that $\supp Z \subseteq R$.
Summing (\ref{tmp:22.05.2012_1}), we obtain
$$
    \| Z \|_1 =  a^2-2t-a = n+1-\Delta = a-\Delta := z \,.
$$
Using this, Gauss sums bound of Lemma \ref{l:K_chi}, inequality $\D\ge 3$, and identity (\ref{tmp:22.05.2012_1}), we get
\begin{equation}\label{est:A_Fourier_Z}
    |\FF{A} (x)|^2 \le 3n+2 < 3a\,, \quad \forall x \neq 0 \,.
\end{equation}
Further,
multiplying (\ref{tmp:22.05.2012_1}) by $(A*A) (x)$ and summing, we have
$$
    \E(A) = \sum_x (2R(x) + (2\cdot A) (x)) (A*A) (x) + \sum_x Z(x) (A*A) (x)
        =
$$
$$
        =
        4t + a + 2z + 4 \sum_x R(x) (2\cdot A) (x) + \sum_x (2\cdot A) (x) Z(x) + \sum_x Z(x) (A*A) (x)
    \,.
$$
Applying Fourier transform and
estimates
(\ref{f:d_zero}), (\ref{f:d_nonzero}), (\ref{f:2A&R}),
 (\ref{est:A_Fourier_Z}), we get
$$
    \E(A) = 2p-2 + a + 2z + 4 \langle 2\cdot A, R \rangle + \langle 2\cdot A, Z \rangle +
        \theta_1 ( \min\{z a, \frac{a^2 z}{p} + a\sqrt{3z} \})
            \le
$$
\begin{equation}\label{f:A_energy}
            \le
                2a^2 - a + \theta_2 (4(a-\Delta)+2a  + \min\{a z, a\sqrt{3z} \}) \,,
\end{equation}
where $|\theta_1|, |\theta_2| \le 1$.
Thus
$\mathcal{E}_1 \le 4z+2a+ \min\{a z, a\sqrt{3z} \}$.
Squaring (\ref{tmp:22.05.2012_1}), summing and using
(\ref{f:d_nonzero}), (\ref{f:2A&R}), we obtain
$$
    \E(A) = 4 t + a + 4z + 2 \langle 2\cdot A, Z \rangle + 4\langle 2\cdot A, R \rangle + \| Z \|_2^2
        = 2a^2 - a + \theta_3 (6(a-\Delta)+2a+z^2) \,,
$$
where $0\le \theta_3 \le 1$.
Hence $\mathcal{E}_1 \le 6z+2a+z^2$.

Finally
$$
    (a-1)^2 + \sum_{x} \mathcal{E}^2_2 (x)
        =
    (a-1)^2 + \sum_{x\neq 0} \left( (A\c A)(x) - 1 \right)^2
                =
            \sum_x \left( (A\c A)(x) - 1 \right)^2
                =
$$
\begin{equation}\label{tmp:28.02.2013_1}
                =
                    \E(A) - 2a^2 + p \,.
\end{equation}
Because of $p=(a-1)^2+\Delta$, the last identity and formula (\ref{f:A_energy_close}),
we obtain $\| \mathcal{E}_2 \|^2_2 = \mathcal{E}_1 + \Delta - a$.

${\bf 3)~}$
Let $a$ be an even number.
Recalling (\ref{f:eps_not_A}) we see that $\eps(x)$ is even for $x\notin -A$.
Using Lemma \ref{l:div_c} with $c=-2$, we get
$$
    \| \eps \|_2^2 = pa-a^2 -a^3 - 2(a-1)d + \o
        \ge
            2(a^2+d) - (a-\o) \,.
$$
After some computations, we obtain
\begin{equation}\label{tmp:16.05.2013_1}
    2 + \D  -2d \ge 5a \,.
\end{equation}
Thus, either $\D=a=n+1$ and $d=-2a+1, -2a$
or $\D = a-1,a-2$ and $d=-2a$.
In the first case $A$ is a perfect difference set.
Now suppose that $A$ is not a perfect difference set and consider the second situation.
Because of $d=-2a$, we
see that
$0\notin A$.
Returning to (\ref{f:C_n^2_fake}), we get $\D=a-2=n-1$.
It means, in particular, that inequality (\ref{tmp:16.05.2013_1}) is actually, equality.
Using formula (\ref{f:div_c'}) of Lemma \ref{l:div_c}, we see that $\eps(x) = 0$ on $A$,
$\eps(x)=-1$ on $-A$ and $\eps(x) = 0,-2$ on the rest.
But $2\cdot A \subseteq N$, hence $A\subseteq R$ or $A\subseteq N$
and thus by (\ref{f:eps_not_A}), we have $\eps(0) = \pm a$ with contradiction.
In particular, we have considered the remain cases $a=4,6$.

Now we should deal with the
situation
$a$ is an odd number
and $0\in A$.
Whence $d=-2a+1$ and we can suppose that $\D \le a-1$ because $A$ is a perfect difference set otherwise.
Using formula (\ref{f:div_c'}) of Lemma \ref{f:div_c} with $c=-2$,
applying the arguments from the proof of Theorem \ref{t:Lev_q}, we obtain
$$
    \| \eps \|_2^2 = pa -a^3 + 3a^2 - 6a + 3
        \ge
            2(a^2-2a+1) + \sum_{k} (k^2+2k) p_k
                \ge
$$
$$
    \ge 2(a^2+2a-1) - p_{-1} + 3 \sum_{k \neq 0,-1,-2} p_k
        =
            2(a^2-2a+1) + 3p - 4 p_{-1} - 3(p_0+p_{-2}) \,.
$$
Here, again $p_k := |\{ x \in \F_p ~:~ \eps(x) = k \}|$, $k\in \Z$.
Clearly, $p_0+p_{-2} \le 2a-1$.
As above, $p_{-1} \le (p-1)/2$.
Hence
$$
    pa -a^3 + 3a^2 - 6a + 3 \ge 2(a^2-2a+1) + p + 2 - 6a + 3 \,.
$$
In other words
\begin{equation*}\label{tmp:17.05.2013}
    \D(a-1) - 2a^2 + 7a - 5 \ge 0
\end{equation*}
with contradiction for $a\ge 5$ because of $A$ is not a perfect difference set by our assumption
and, hence, $\D \le a-1$.
This completes the proof.
$\hfill\Box$

\bigskip

Clearly, the third part of Theorem \ref{t:A+A=R} combining with Proposition \ref{pr:pd_sets}
give the second part of the main Theorem \ref{t:main_intr}.

\bigskip

\Note
{
The method of the proving Theorem \ref{t:A+A=R} above is analytical, thus we can say
something about the structure of the set $A$ in a little bit more general situation.
It was the reason to include the second part in the theorem.
More precisely, suppose that a formula
    \begin{equation}\label{f:gen_Z}
        (A*A) (x) = 2R(x) + (2\cdot A)(x) + Z(x)
    \end{equation}
with "small"\, function $Z$ takes place.
Then
summing (\ref{f:gen_Z}), we have
$$
    p-1 + |A| - \| Z\|_1 \le |A|^2 \le p-1 + |A| + \| Z\|_1 \,,
    \quad \mbox{ and } \quad \langle Z \rangle = |A|^2 - 2t-a \,.
$$
Further, for all $x\neq 0$ the following holds
$$
    |\FF{A} (x)|^2 \le \sqrt{p}+1+|A| + \| Z\|_1 \,,
$$
and
$$
    \E(A) = \sum_x (2R(x) + (2\cdot A) (x)) (A*A) (x) + \sum_x Z(x) (A*A) (x)
        \le
$$
$$
        \le
        4t + 5a + 3\| Z\|_1 + \sum_x Z(x) (A*A) (x)
            \le
$$
\begin{equation}\label{tmp:28.02.2013_2}
            \le
                4t + 6a + 3\| Z\|_1 + (\sqrt{p}+1+|A| + \| Z\|_1)^{1/2} |A|^{1/2} \| Z \|_2 \,,
\end{equation}
provided by $\langle Z \rangle \le pa^{-1}$, say.
Similarly, if we put
\begin{equation*}\label{f:A_difference_close}
        (A\c A) (x) = (|A|-1) \d_0 (x) + 1 + \mathcal{E}_2 (x) \,,
\end{equation*}
Then, as in (\ref{tmp:28.02.2013_1}), we get by (\ref{tmp:28.02.2013_2}) that
$$
    \sum_{x} \mathcal{E}^2_2 (x)
                =
            \sum_x \left( (A\c A)(x) - 1 \right)^2 - (a-1)^2
                =
                    \E(A) - 2a^2 + p - (a-1)^2
                        \le
$$
$$
    \le
        3p-3 + 8a - 3a^2 + 3\| Z\|_1 + (\sqrt{p}+1+|A| + \| Z\|_1)^{1/2} |A|^{1/2} \| Z \|_2 \,.
$$
Thus, if $|R \Delta (A\hat{+} A)|$ is small than $A$ is structured and close to a perfect difference set.
}

\refstepcounter{section}
\label{sec:appendix}

\begin{center}
{\bf \large
    \arabic{section}.
    Concluding remarks}
\end{center}

Now we describe an alternative way to obtain the results of sections \ref{sec:proof}, \ref{sec:proof2}.

Our aim is
to define an analog of "Fourier transform" relatively to a multiplicative character $\chi$.
Put
$$
    \mathrm{f} (x) = \frac{1}{\sqrt{p}} \left( \chi(x) - \frac{1}{\sqrt{p}} \right) \,.
$$
Clearly, $\mathrm{f}$ is a perfect difference function.
Note that
\begin{equation}\label{f:uniform_bound_f}
    \| \mathrm{f} \|_{\infty} \le \frac{1}{\sqrt{p}} \left( 1+ \frac{1}{\sqrt{p}} \right) \,.
\end{equation}
Let $g$ be an arbitrary function.
We write $g_s(x)$ for $g(x+s)$.
By Lemma \ref{l:K_chi}, we have
\begin{equation}\label{f:onb_f}
    \langle \mathrm{f}_s, \mathrm{f}_t \rangle = \d_{s,t}
        \quad \mbox{ and } \quad
            \sum_y \mathrm{f}_y (s) \ov{\mathrm{f}_y (t)} = \d_{s,t} \,,
\end{equation}
where $s,t \in \F_p$ are arbitrary.
Further $\langle \mathrm{f}_s, 1 \rangle = \langle 1, \mathrm{f}_s \rangle = -1$ and
\begin{equation}\label{}
    (\mathrm{f}_s \c \ov{\mathrm{f}}_t) (x) = \sum_z \mathrm{f}_s (z) \ov{\mathrm{f}_t (x+z)}
        = \langle \mathrm{f}_s, \mathrm{f}_{t+x} \rangle = \d_{s,t+x} \,.
\end{equation}

For an arbitrary function $g$ write $g^{\la} (x) = g(\la x)$.
We have
\begin{equation}\label{}
    \mathrm{f}^{\la}_s (x) = \chi(\la) \mathrm{f}_{s/\la} (x) + \frac{\chi(\la)-1}{p} \,,
\end{equation}
where $\la$ is any nonzero element.
Applying formula (\ref{f:K_chi_convolution}) once more, we obtain
\begin{equation}\label{f:*_f}
    (\mathrm{f}^{\la}_s \c \ov{\mathrm{f}}^{\mu}_t ) (x) = \chi(\la) \ov{\chi(\mu)} \cdot \d_{s/\la,x+t/\mu} + \frac{1-\chi(\la) \ov{\chi(\mu)}}{p} \,.
\end{equation}

Now we can define the "Fourier transform"\, corresponding to the character $\chi$.

\Def
{
Let $\_phi : \F_p \to \C$ be a function.
By $\t{\_phi}$ denote the function
\begin{equation}\label{def:tilde_transform}
    \t{\_phi} (x) = \sum_y \_phi(y) \ov{\mathrm{f}_x (y)} = \langle \_phi, \mathrm{f}_x \rangle = (\_phi \c \ov{\mathrm{f}}) (x) \,.
\end{equation}
}

For example $\t{1} = -1$, $\t{\d}_s (x) = \ov{\mathrm{f}_s (x)}$ and $\t{\mathrm{f}}_s (x) = \d_s (x)$.

Because of $\mathcal{D}$ is closed under convolutions $*$, $\c$, and $\FF{\mathrm{f}} (x) \neq 0$, $x\in \F_p$,
we see that $\_phi \in \mathcal{D}$ iff $\t{\_phi} \in \mathcal{D}$
(formula (\ref{f:tilde_2}) below
also implies the fact).
The
next
lemma follows from the definitions.

\Lemma
{\it
    Let $\_phi, \psi : \F_p \to \C$ be any functions.
    Then
\begin{equation}\label{f:tilde_0}
    (\_phi_s)\, \t{} = (\t{\_phi})_{-s} \,, \quad \forall s\in \F_p
        \,.
\end{equation}
\begin{equation}\label{f:tilde_1}
    (\ov{\t{\_phi}}) \, \t{} \, (x) = \ov{\_phi} (x)
        \quad \quad \mbox{ and } \quad \quad
                (\t{\_phi}) \, \t{} \, (x) = \_phi (x)
            \quad \mbox{ if } \quad \chi \in \R \,.
\end{equation}
\begin{equation}\label{f:tilde_1.5}
    (\_phi^{\la} )\, \t{} \, (x) = \chi(\la) \t{\_phi}^{\la} (x) + \frac{\langle \_phi \rangle (\chi(\la)-1)}{p}
        \,, \quad \quad \forall \la \neq 0 \,.
\end{equation}
\begin{equation}\label{f:tilde_2}
    \_phi (x) = \sum_y \t{\_phi} (y) \mathrm{f}_x (y) \,.
\end{equation}
\begin{equation}\label{f:tilde_3}
    \sum_x \_phi (x) \ov{\psi (x)} = \sum_y \t{\_phi} (y) \ov{\t{\psi} (y)}
\end{equation}
\begin{equation}\label{f:tilde_5}
    (\_phi \c \psi)\, \t{} \, (x) = (\_phi * \t{\psi}) (x) \,.
\end{equation}
Further for any $\la,\mu\neq 0$, we have
\begin{equation}\label{f:tilde_4*}
    (\_phi^{\la} \c \ov{\psi^{\mu}}) (x) =
            \chi(\la) \ov{\chi(\mu)} \cdot (\ov{\t{\psi}^{\mu}} \c \t{\_phi}^{\la}) (x)
                +
            \frac{1-\chi(\la) \ov{\chi(\mu)}}{p} \cdot \langle \_phi \rangle \langle \ov{\psi} \rangle
        \,.
\end{equation}
In particular
\begin{equation}\label{f:tilde_4}
    (\_phi \c \ov{\psi}) (x) = (\t{\_phi} \c \ov{\t{\psi}}) (-x) = (\ov{\t{\psi}} \c \t{\_phi}) (x)
\end{equation}
and
\begin{equation}\label{f:tilde_4'}
    (\_phi * \ov{\psi} ) (x) = \chi(-1) (\t{\_phi} * \ov{\t{\psi}}) (-x)
        + \frac{1-\chi(-1)}{p} \langle \_phi \rangle \langle \ov{\psi} \rangle \,.
\end{equation}
}
\label{l:tilde_Fourier}
\Proof
We just need to check (\ref{f:tilde_4*})
because another formulas
are almost trivial.
Applying (\ref{f:tilde_1.5}), (\ref{f:tilde_2}) and the fact $\langle \mathrm{f}_s \rangle = -1$, $s\in \F_p$, we obtain
$$
    (\_phi^{\la} \c \ov{\psi^{\mu}}) (x) = \sum_z \_phi^{\la} (z) \ov{\psi^{\mu} (z+x)}
        =
$$
$$
        =
            \sum_{y,y',z} \left( \chi(\la) \t{\_phi}^{\la} (y) + \frac{\langle \_phi \rangle (\chi(\la)-1)}{p} \right)
                \cdot
                     \ov{\left( \chi(\mu) \t{\psi}^{\mu} (y') + \frac{\langle \psi \rangle (\chi(\mu)-1)}{p} \right)}
                        \mathrm{f}_z (y) \ov{\mathrm{f}_{z+x} (y')}
                            =
$$
$$
    =
        \chi(\la) \ov{\chi(\mu)} \cdot \sum_y \t{\_phi}^{\la} (y) \ov{\t{\psi}^{\mu} (y-x)}
            +
                \frac{1-\chi(\la) \ov{\chi(\mu)}}{p} \cdot \langle \_phi \rangle \langle \ov{\psi} \rangle
                    =
$$
$$
    =
         \chi(\la) \ov{\chi(\mu)} \cdot (\ov{\t{\psi}^{\mu}} \c \t{\_phi}^{\la}) (x)
                +
            \frac{1-\chi(\la) \ov{\chi(\mu)}}{p} \cdot \langle \_phi \rangle \langle \ov{\psi} \rangle \,.
$$
This concludes the proof.
$\hfill\Box$

In particular, formula (\ref{f:tilde_4}) says that $\E(A) = \E (\t{A})$.
The transform above can be used to obtain the results of sections \ref{sec:proof}, \ref{sec:proof2}.
Indeed, the arguments here are just calculation of $\t{A}$, $\t{B}$.
Similarly, one can define the "Fourier transform"\, relatively any perfect difference set
as well as a function from $\mathcal{D}$.

\bigskip

There is a general question about "Fourier coefficients"\, of subsets of $\F_p$.
Partial answer to the question give us the proof of the statements of sections \ref{sec:proof}, \ref{sec:proof2}.

Suppose that $A\subset \F_p$ is a set,
and
$\chi$ is an arbitrary character.
Consider the function
$$
    E(x) = E_A (x) = E_{A,\chi} (x):= (A\c \chi) (x) \,.
$$
Alternatively, one can take convolutions of $A$ with the function $\mathrm{f}$ above, that is $\t{A}$.
Our question is the following : what can we say nontrivial about the function $E(x)$?
There is a list of simple properties of $E(x)$.
Clearly,
$$
    \langle E_A \rangle = 0 \quad \mbox{ and } \quad \| E_A \|_2^2 = p|A|-|A|^2 \,.
$$
Higher moments of $E$ can be
estimated
using Weil's
bound
(\ref{f:Weil}) and Lemma \ref{l:E_k-identity}
as it was done in the proof of Corollary \ref{cor:Sarkozy}.
Of course, there is the "cocycle"\, property of $E(x)$, namely,
$$
    E_A (xy) = \chi(x) E_{x^{-1} A} (y) \,, \quad \quad x\neq 0 \,.
$$
In particular, it means that $|E_A (x)|$ is
fixed
by the multiplier group of $A$, that is the set $\{ x\in \F_p ~:~ xA=A \}$.
Further,
$$
    \sum_{x,y} \chi(x+y) A(x) E_A (y) = \sum_z (A\c \chi)^2 (z) = p|A|-|A|^2
$$
and hence
\begin{equation}\label{f:l_opt}
    |\sum_{x,y} \chi(x+y) A(x) E_A (y)| \gg \sqrt{p} \| A \|_2 \| E_A \|_2
\end{equation}
provided by $|A| \ll p$.
In particular, formula (\ref{f:l_opt}) shows that inequality (\ref{f:rho_of_sum_1})
of Lemma \ref{l:rho_of_sum} is tight.
Using estimate (\ref{f:uniform_bound_f}) and Lemma \ref{l:tilde_Fourier} one can easily prove
an uncertainty principle for any function $g$ (in particular for $A$), namely
$$
    p \left( 1 + \frac{1}{\sqrt{p}} \right)^{-2} \le |\supp g| \cdot |\supp \t{g}| \,.
$$
If $A$ is a Singer's perfect difference set then the character $\chi$
can be represented as the convolution $(A*E_{A,\chi})(x)$.
Indeed
$$
    (A*E_A) (x) = (A*(A\c \chi)) (x) = (\chi * (A\c A)) (x) = (\chi * ((|A|-1) \d_0 + 1)) (x)
        =
            (|A|-1) \chi (x)\,.
$$
Finally, the size of each level set
$$
    L_c := \{ x \in \F_p ~:~ E_A (x) = c \}\,, \quad \quad c\in \F_p
$$
is bounded by $(p-1)/2$ because of $E (x) \pmod p$ is a nonzero polynomial from $\F_p [x]$ of degree
$(p-1)/2$
with the leading term $|A| x^{(p-1)/2}$.
Are there further properties of $E(x)$?
For example, what can we say about the maximum/minimum value of $E(x)$?

\bigskip

\no{Division of Algebra and Number Theory,\\ Steklov Mathematical
Institute,\\
ul. Gubkina, 8, Moscow, Russia, 119991\\}
and
\\
Delone Laboratory of Discrete and Computational Geometry,\\
Yaroslavl State University,\\
Sovetskaya str. 14, Yaroslavl, Russia, 150000\\
and
\\
IITP RAS,  \\
Bolshoy Karetny per. 19, Moscow, Russia, 127994\\
{\tt ilya.shkredov@gmail.com}

\end{document}